\documentclass[11pt,reqno]{amsart}
\usepackage{color}

\newtheorem{theorem}{Theorem}
\newtheorem{definition}{Definition}
\newtheorem{proposition}[theorem]{Proposition}

\newtheorem{remark}{Remark}
\newtheorem{lemma}{Lemma}

\newfont{\bb}{msbm10 at 12pt}
\def\pf{{\textit {Proof :} }}

\def\R{\hbox{\bb R}}

\def\SE{{\mathbb{S} \E}}

\def\Spin{{\mathop{\rm Spin}}}  
\def\bS{{\mathbb S}}

\def\E{\mathcal{E}}

\def\<{\langle}     
\def\>{\rangle}     
\def\Tilde{\widetilde}
\def\div{{\rm div}}

\def\TL{\Tilde{\mathcal{L}}}
\def\L{\mathcal{L}}

\newcommand{\bal}{\begin{align}}      \newcommand{\eal}{\end{align}}
\newcommand{\ba}{\begin{array}}      \newcommand{\ea}{\end{array}}
\newcommand{\bc}{\begin{center}}     \newcommand{\ec}{\end{center}}
\newcommand{\be}{\begin{enumerate}}  \newcommand{\ee}{\end{enumerate}}
\newcommand{\beQ}{\begin{eqnarray*}} \newcommand{\eeQ}{\end{eqnarray*}}
\newcommand{\bi}{\begin{itemize}}    \newcommand{\ei}{\end{itemize}}
\newcommand{\bt}{\begin{tabular}}    \newcommand{\et}{\end{tabular}}
\newcommand{\bdm}{\begin{displaymath}} \newcommand{\edm}{\end{displaymath}}

    \newcommand{\sm}{\bS\!\!\!/\,\!}
\newcommand{\D}{D\!\!\!\!/\,}
\newcommand{\EMD}{\mathcal{D}\!\!\!\!/\,}
\newcommand{\nb}{\nabla\!\!\!\!/\,}
\newcommand{\mult}{\gamma\!\!\!/}

\def\qed{\hfill{q.e.d.}\smallskip\smallskip}

\begin{document}

\title[An Alexandrov theorem in Minkowski spacetime]{An Alexandrov theorem in Minkowski spacetime}

\author{Oussama Hijazi}
\address[Oussama Hijazi]{Institut {\'E}lie Cartan,
Universit{\'e} de Lorraine, Nancy,
B.P. 239,
54506 Vand\oe uvre-L{\`e}s-Nancy Cedex, France.}
\email{Oussama.Hijazi@univ-lorraine.fr}

\author{Sebasti{\'a}n Montiel}
\address[Sebasti{\'a}n  Montiel]{Departamento de Geometr{\'\i}a y Topolog{\'\i}a,
Universidad de Granada,
18071 Granada,  Spain.}
\email{smontiel@goliat.ugr.es}

\author{Simon Raulot}
\address[Simon Raulot]{Laboratoire de Math\'ematiques R. Salem
UMR $6085$ CNRS-Universit\'e de Rouen
Avenue de l'Universit\'e, BP.$12$
Technop\^ole du Madrillet
$76801$ Saint-\'Etienne-du-Rouvray, France.}
\email{simon.raulot@univ-rouen.fr}

\begin{abstract}
In this paper, we generalize a theorem {\it \`a la Alexandrov} of Wang, Wang and Zhang \cite{WWZ} for closed codimension-two spacelike submanifolds in the Minkowski spacetime for an adapted CMC condition. 
\end{abstract}



\thanks{The second author was partially  
supported by a Spanish MEC-FEDER grant No. MTM2011-22547}

\date{\today}   

\maketitle 
\pagenumbering{arabic}
 

\section{Introduction}


The classical Alexandrov theorem \cite{A} asserts that the only compact embedded hypersurfaces with constant mean curvature (CMC) in the Euclidean space are the round spheres. Natural generalizations of this result has been obtained for such hypersurfaces in the hyperbolic space and the open hemisphere \cite{MR} as well as in some warped product manifolds (see \cite{Mo} or \cite{Br} more recently). 

Since Euclidean space and hyperbolic space arise as spacelike hypersurfaces in the Minkowski spacetime, a natural question is whether one can obtain these two results as a particular case of a more general result concerning codimension-two submanifolds in the Minkowski spacetime. From the point of view of submanifolds theory, a natural analogue of the CMC condition for higher codimensional submanifolds is the parallel mean curvature condition. From the general relativity point of view, the most relevant physical phenomenon is the divergence of light rays emanating from a codimension-two submanifolds. More precisely, the causal future or past of a geometric object is of great importance. It is interesting to characterize when a surface lies in the null hypersurface generated by a ``round sphere''. These are called ``shearfree'' null hypersurfaces in general relativity literature, and are analogues of umbilical hypersurfaces in Riemannian geometry. This type of condition has recently been introduced by Wang, Wang and Zhang \cite{WWZ} and is described as follow. 

Consider $\Sigma^n$ a codimension-two spacelike orientable submanifold in a $(n+2)$-dimensional Lorentzian manifold $(\E^{n+1,1},\<\;,\;\>)$. We will represent by ${\mathcal H}$ the mean curvature vector field on $\Sigma^n$, defined as 
$$
{\mathcal H}={\rm tr\,}{\rm II},
$$ 
where ${\rm II}$ is the second fundamental form of the immersion $\Sigma^n\hookrightarrow\E^{n+1,1}$ given by
\begin{eqnarray*}
{\rm II}(X,Y)=\Tilde{\nabla}_XY-\nabla^\Sigma_XY
\end{eqnarray*}
for all $X,Y\in\Gamma(T\Sigma)$ and where $\widetilde{\nabla}$ (resp. $\nabla^\Sigma$) is the Levi-Civita connection of $\mathcal{E}^{n+1,1}$ (resp. $\Sigma^n$). Since the normal space at each point of $\Sigma^n$ is a Lorentzian
plane, it can be spanned by two future-directed null normal vector field ${\mathcal L}_+$ and ${\mathcal L}_-$ normalized in such a way that $\langle{\mathcal L}_+,{\mathcal L}_-\rangle=-2$. In this frame, the second fundamental form can be expressed as
\begin{eqnarray*}
{\rm II}(X,Y)=\frac{1}{2}\chi_+(X,Y)\mathcal{L}_-+\frac{1}{2}\chi_-(X,Y)\mathcal{L}_+
\end{eqnarray*}
where $\chi_\pm$ is the null second fundamental form with respect to ${\mathcal L}_{\pm}$ defined by 
\begin{eqnarray}\label{NSFF}
\chi_\pm(X,Y)=\<\widetilde{\nabla}_X\mathcal{L}_\pm,Y\>
\end{eqnarray}
for all $X$, $Y\in\Gamma(T\Sigma)$. 
We denote by $\theta_\pm={\rm tr\,}\chi_\pm$ the so-called future-directed null expansions of ${\mathcal H}$ which measure the area growth when $\Sigma^n$ varies in the corresponding directions. It is then clear that
\begin{eqnarray*}
\mathcal{H}=\frac{1}{2}\theta_+{\mathcal L}_-+\frac{1}{2}\theta_-{\mathcal L}_+\quad\text{and}\quad|\mathcal{H}|^2= -\theta_+\theta_-.
\end{eqnarray*}
If $\theta_+$ and $\theta_-$ are both negative, the submanifold will be called a {\em trapped} submanifold. A codimension-two submanifold with $\theta_+=0$ or $\theta_-=0$ is called  a {\em marginally trapped} submanifold. Remark that if $\Sigma^n$ is trapped or marginally trapped, then the mean curvature vector ${\mathcal H}$ is a causal vector at each point. This is why that if the mean curvature field ${\mathcal H}$ is spacelike everywhere, $\Sigma$ will be refer to as an {\em untrapped} submanifold.  

In the case where $\Sigma^n$ spans a spacelike hypersurface in the Lorentzian manifold, that is, when there exists a spacelike hypersurface $\Omega^{n+1}$ immersed in $\E^{n+1,1}$ such that $\partial\Omega^{n+1}=\Sigma^n$,  the normal null vector fields ${\mathcal L}_+$ and ${\mathcal L}_-$ may be ordered in such a way that they project onto directions tangent to $\Omega^{n+1}$ which are respectively {\em outer} and {\em inner} normal at each point of $\Sigma^n$. In other words, if $N$ is an inner normal unit vector field on $\Sigma^n$ tangent to $\Omega^{n+1}$ and $T$ is a future-directed timelike normal to $\Omega^{n+1}$ in $\E^{n+1,1}$, we put
$$
{\mathcal L}_+=T-N,\qquad {\mathcal L}_-=T+N.
$$
The second fundamental form of $\Sigma^n$ in $\E^{n+1,1}$ is given in terms of the Lorentzian basis of the normal bundle provided by the hypersurface $\Omega^{n+1}$ by
\begin{eqnarray*}
{\rm II}(X,Y)=\langle AX,Y\rangle N+\langle BX,Y\rangle T,
\end{eqnarray*}
for all $X,Y\in\Gamma(T\Sigma)$ and where 
\begin{eqnarray}\label{shape1}
\<BX,Y\>=\<\Tilde{\nabla}_X T,Y\>
\end{eqnarray}
and 
\begin{eqnarray}\label{shape2}
AX:=-\nabla^\Omega_XN
\end{eqnarray}
denote respectively the shape operators of $\Omega^{n+1}$ in $\mathcal{E}^{n+1,1}$ and $\Sigma^n$ in $\Omega^{n+1}$. Here $\nabla^\Omega$ denotes the Levi-Civita connection of the Riemannian metric $\langle\;,\;\rangle$ on $\Omega$. The mean curvature vector field $\mathcal{H}$ of $\Sigma$ in $\E$ can be re-expressed by:
\begin{eqnarray*}
\mathcal{H}=H N+K T,
\end{eqnarray*}
where $H={\rm tr\,}A$ is the mean curvature of $\Sigma^n$ in $\Omega^{n+1}$ and $K$ is the trace on $\Sigma^n$ of the shape operator $B$ of $\Omega^{n+1}$ in ${\mathcal E}^{n+1,1}$, that is $K={\rm tr}_\Sigma\,B$. The norm of ${\mathcal H}$ can also be re-expressed as
\begin{eqnarray*}
|\mathcal{H}|^2= H^2-K^2,
\end{eqnarray*}
with $\theta_\pm=K\pm H$ are the future-directed null expansions of ${\mathcal H}$. 
The spacelike codimension-two submanifolds with $\theta_+<0$ (respectively, $\theta_-<0$) are referred to as {\em outer} (respectively, {\em inner}) {\em trapped} submanifolds. For these reasons, a codimension-two untrapped submanifold which bounds a compact connected spacelike hypersurface $\Omega$ in $\E$ and which is mean convex in $\Omega$ will be referred to as an {\it outer untrapped }submanifold. It is worth noting that round spheres in Euclidean slices of the Minkowski spacetime are outer untrapped submanifolds as well as large radial spheres in asymptotically flat spacelike hypersurfaces.

Now recall that a closed oriented spacelike codimension-two submanifold $\Sigma$ in a $(n+2)$-dimensional Lorentzian manifold is said to be {\it torsion-free} with respect to a null normal vector field $\mathcal{L}$ along $\Sigma$ if the connection one-form $\zeta_\mathcal{L}$ defined by 
\begin{eqnarray*}
\zeta_\mathcal{L}(X)=\frac{1}{2}\<\widetilde{\nabla}_X\mathcal{L},\underline{\mathcal{L}}\>
\end{eqnarray*}
is zero for all $X\in\Gamma(T\Sigma)$. Here $\underline{\mathcal{L}}$ is another null normal such that $\<\mathcal{L},\underline{\mathcal{L}}\>=-2$. This condition is easily seen to be equivalent to $(\widetilde{\nabla}\mathcal{L})^\perp=0$ where $(\,.\,)^\perp$ denotes the normal component.
\begin{definition}
A codimension-two submanifold of a Lorentzian manifold is said to have {\it constant normalized null curvature} (CNNC) if there exists a future null normal vector field $\mathcal{L}$ such that $\Sigma$ is torsion-free with respect to $\L$ and $\<\mathcal{H},\mathcal{L}\>$ is a constant. 
\end{definition}
 
Obviously, the CNNC assumption reduces to the CMC assumption when $\Sigma$ lies in a totally geodesic spacelike hypersurface of a Lorentzian manifold. In \cite{WWZ}, the authors prove among other things the following result 
\begin{theorem}\cite{WWZ}
Let $\Sigma^n$ be a future incoming null smooth, closed, embedded, spacelike codimension-two submanifold in the Minkowski spacetime. Suppose that $\Sigma$ has CNNC with respect to a future incoming null normal vector field $\mathcal{L}$ and $\<\mathcal{H},\mathcal{L}\>>0$. Then $\Sigma$ lies in a shearfree null hypersurface. 
\end{theorem}
 
A closed, spacelike codimension-two submanifold $\Sigma$ in a static spacetime is future (resp. past) incoming null smooth if the future (resp. past) incoming null hypersurface of $\Sigma$ intersects a totally geodesic time-slice of the spacetime at a smooth orientable hypersurface. Moreover a null hypersurface $\mathcal{C}$ in $\mathcal{E}$ is shearfree if there exists a spacelike hypersurface $\mathcal{S}\subset\mathcal{C}$ (so a codimension-two submanifold of $\mathcal{E}$) such that the null second fundamental form $\chi$ of $\mathcal{S}$ with respect to some null normal $\mathcal{L}$ satisfies $\chi=f \<\;,\;\>$ for function f and where $\<\;,\;\>$ is the Riemannian metric induced on $\mathcal{S}$ by the Lorentzian one. Note that being shearfree is a property of the null hypersurface. 

Their proof relies on Heintze-Karcher-type inequalities and on a certain monotonicity formula which strongly use the incoming null smoothness of $\Sigma$ in the spacetime.

In this paper, we will use the spinorial approach developed by the first two authors and X. Zhang \cite{HMZ} (see also \cite{HMRo}) to generalize this result by relaxing the assumption on the incoming null smoothness of $\Sigma$. More precisely, we have 
\begin{theorem}\label{AlexMink}
Let $\Sigma^n$ be an untrapped codimension-two submanifold in the Minkowski spacetime. Suppose that $\Sigma$ has CNNC with respect to a future null normal vector field $L$. Then $\Sigma$ lies in a shearfree null hypersurface. 
\end{theorem}

The proof of this result relies on a conformal eigenvalue estimate for a Dirac-type operator acting on spinors of $\Sigma$ and especially on a careful treatment of its equality case. Note that our method does not seem to apply directly in the general context treated in \cite{WWZ} since we need the existence of a particular section of the spinor bundle, namely an extrinsic imaginary Codazzi spinor (see (\ref{EBK})). 

\begin{remark}{\rm
It should be point out that in \cite{HMR2} we investigate the rigidity of time flat submanifolds in Minkowski spacetime using a similar method. However the Dirac-type operator we used in this paper is not the one study in the present work. }
\end{remark}


\section{Preliminaries on spacetime geometry}



\subsection{The Einstein equation, the dominant energy condition and geometry of codimension-two spacelike submanifolds} 


Let $(\E^{n+1,1},\<\;,\;\>)$ be a time-oriented spacetime satisfying the Einstein field equations, that is $\E^{n+1,1}$ is an $(n+2)$-dimensional smooth manifold endowed with a smooth Lorentzian metric $\<\;,\;\>$ such that
\begin{eqnarray}\label{Einstein}
\Tilde{Ric}-\frac{1}{2}\Tilde{R}\<\;,\;\>=\mathcal{T},
\end{eqnarray}
where $\Tilde{R}$ (respectively, $\Tilde{Ric}$) denotes the scalar curvature (respectively, the Ricci curvature) of ($\E^{n+1,1},\<\;,\;\>)$ and $\mathcal{T}$ is the energy-momentum tensor which describes the matter content of the ambient spacetime. 

Let $M^{n+1}$ be an immersed spacelike hypersurface of $\E^{n+1,1}$ whose induced Riemannian metric is also denoted by $\langle\;,\;\rangle$. Let $T$ be the future-directed timelike unit vector field normal to $M$ and denote by $B$ the associated shape operator given by (\ref{shape1}). Then the Gau{\ss}, Codazzi and Einstein equations provide {\it constraint equations} on $M^{n+1}$ 
\begin{eqnarray}
\mu & = & \frac{1}{2}\big(R-|B|^2+({\rm tr\,}B)^2\big) \label{constraint1}\\   
J & = & \div\big(B-({\rm tr\,}B)I\big) \label{constraint2}
\end{eqnarray}
where $R$ is the scalar curvature of $(M^{n+1},\langle\;,\;\rangle)$, $|B|^2$ and ${\rm tr\,}B$ denote the squared norm and the trace of $B$ on $M^{n+1}$ with respect to $\langle\;,\;\rangle$. The quantity $\mu$ (respectively, $J$) is the energy (respectively, the momentum) density of the matter fields given by 
\begin{eqnarray*}
\mu = \mathcal{T}(T,T)\quad\big(\text{respectively,}\quad \langle J,v\rangle = \mathcal{T}(T,v)
\end{eqnarray*} 
for each spacelike vector $v$ tangent to $M^{n+1}\big)$. A triplet $(M^{n+1},\langle\;,\;\rangle,B)$ which satisfies the constraint equations (\ref{constraint1}) and (\ref{constraint2}) for given densities $\mu$ and $J$ is called an {\it initial data set} for the Einstein field equations.

Let us suppose from now on that the spacetime satisfies the dominant energy condition (DEC), that is, its energy-momentum tensor $\mathcal{T}$ has the property that the vector field dual to the one-form $-\mathcal{T}(u,\cdot)$ is a future-directed causal vector of $T\E$, for every future-directed causal vector $u\in\Gamma(T\E)$. This implies in particular that the following inequality holds
\begin{equation}\label{DEC}\tag{${\rm DEC}$}
\mu |v|\geq |\langle J,v\rangle|,
\end{equation}
for all $v\in TM$.


\subsection{Spin geometry of spacelike hypersurfaces in $\E^{n+1,1}$} 


From now we assume that the hypersurface $M^{n+1}$ is spin and let $\Omega^n$ be a domain in $M^{n+1}$. In this situation, the spinor bundle over $\E$ exists globally along $M$. Indeed, if $\Spin(\E)$ stands for the $\Spin_{n+1,1}$-bundle of spinorial frame locally defined in a neighborhood of $M$ in $\E$ then the associated complex spinor bundle $\SE$ is defined by:
\begin{eqnarray*}
\SE:=\Spin(\E)\times_{\Tilde{\gamma}_{n+1,1}}\mathcal{S}_{n+2},
\end{eqnarray*}
where $\Tilde{\gamma}_{n+1,1}$ is the complex representation of the group $\Spin_{n+1,1}$ and $\mathcal{S}_{n+2}$ is the $\Spin_{n+1,1}$-module of complex dimension $2^{[\frac{n+2}{2}]}$. On the other hand, the existence of a unit timelike vector $T$ normal to $M$ (and so to $\Omega$) allows to define the restricted spinor bundle $\sm \Omega$ by:
\begin{eqnarray*}
\sm \Omega=\SE_{|\Omega}:=\Spin(\Omega)\times_{\Tilde{\gamma}_{n+1,1}\circ\eta}\mathcal{S}_{n+2}
\end{eqnarray*} 
where $\eta$ is the natural inclusion $\Spin_{n+1}\subset \Spin_{n+1,1}$ and where $\Spin(\Omega):=\Spin(\E)_{|\Omega}$. The natural action of $\omega\in\mathbb{C}l(\E)$, an element of the complex Clifford bundle over $\E$, on a spinor field $\psi\in\Gamma\big(\sm \Omega\big)$ will be denoted by $\Tilde{\gamma}(\omega)\psi$. This action induces a Clifford multiplication on $\sm\Omega$ denoted by $\mult^\Omega$ and related to $\Tilde{\gamma}$ by:
\begin{eqnarray}\label{idmultclif}
\mult^\Omega(X)\psi=i\Tilde{\gamma}(X)\Tilde{\gamma}(T)\psi,
\end{eqnarray}
for all $X\in\Gamma(T\Omega)$ and $\psi\in\Gamma\big(\sm\Omega\big)$. According to \cite{TB-Baum}, the spinor bundle $\sm\Omega$ carries a $\mathrm{Spin}_{n+1,1}$-invariant inner product $(\;,\;)$ such that
\begin{eqnarray*}
(\Tilde{\gamma}(X)\varphi,\psi)=(\varphi,\Tilde{\gamma}(X)\psi), 
\end{eqnarray*}
for all $X\in\Gamma(T\E_{|\Omega})$ and $\varphi,\psi\in\Gamma\big(\sm\Omega\big)$ but which is not positive-definite. However, if one let:
\begin{eqnarray*}
\<\varphi,\psi\>:=(\Tilde{\gamma}(T)\varphi,\psi),
\end{eqnarray*}
it defines a $\Spin_{n+1}$-invariant positive-definite inner product such that
\begin{eqnarray}\label{MultRule}
\<\Tilde{\gamma}(X)\varphi,\psi\>=-\<\varphi,\Tilde{\gamma}(X)\psi\>\qquad\text{and}\qquad\<\Tilde{\gamma}(T)\varphi,\psi\>=\<\varphi,\Tilde{\gamma}(T)\psi\>
\end{eqnarray}
for all $X\in\Gamma(T\Omega)$ and $\varphi,\psi\in\Gamma\big(\sm\Omega\big)$. 

From a Lorentzian point of view, the Gau{\ss} formula gives a relation between the space-time connection $\Tilde{\nabla}$ and the one induced on $T\Omega$ denoted by $\nabla^\Omega$. Namely, we have
\begin{eqnarray*}
\widetilde{\nabla}_X Y=\nabla^\Omega_X Y+\<BX,Y\>T
\end{eqnarray*}
for all $X,Y\in\Gamma(T\Omega)$ and where $B$ is defined by (\ref{shape1}). The spin Gau{\ss} formula gives the counterpart of this formula in the spinorial setting:
\begin{eqnarray}\label{spingaussformule}
\widetilde{\nabla}_X\psi=\nb^\Omega_X\psi+\frac{1}{2}\Tilde{\gamma}\big(BX\big)\Tilde{\gamma}(T)\psi
\end{eqnarray}
for all $X\in\Gamma(T\Omega)$, $\psi\in\Gamma\big(\sm\Omega\big)$ and where $\Tilde{\nabla}$ and $\nb^\Omega$ correspond to the spin Levi-Civita connections obtained by lifting to the spinor bundle $\sm\Omega$ the Lorentzian and Riemannian connections $\Tilde{\nabla}$ and $\nabla^\Omega$. From this identity, we easily compute that 
\begin{eqnarray*}
\nb^\Omega_X\big(\Tilde{\gamma}(T)\psi\big)=\Tilde{\gamma}(T)\nb^\Omega_X\psi
\end{eqnarray*}
and since the connection $\widetilde{\nabla}$ and the inner product $(\;,\;)$ are compatible, we also deduce that
\begin{eqnarray*}
X\<\psi,\varphi\>  = \<\nb^\Omega_X\psi,\varphi\>+\<\psi,\nb^\Omega_X\varphi\>
\end{eqnarray*}
for all $X\in\Gamma(T\Omega)$ and $\psi,\varphi\in\Gamma\big(\sm\Omega\big)$. On the other hand, combining the compatibility of the spin covariant derivative $\widetilde{\nabla}$ with the Clifford product $\widetilde{\gamma}$, that is
\begin{eqnarray*}
\Tilde{\nabla}_X(\Tilde{\gamma}(Y)\psi)=\Tilde{\gamma}(\Tilde{\nabla}_XY)\psi+\Tilde{\gamma}(Y)\Tilde{\nabla}_X\psi
\end{eqnarray*}
for $X,Y\in\Gamma(T\Omega)$ and $\psi\in\Gamma(\sm\Omega)$, with (\ref{spingaussformule}) imply that
\begin{eqnarray}\label{compatibility1}
\nb^\Omega_X(\Tilde{\gamma}(Y)\psi)=\Tilde{\gamma}(\nabla^\Omega_XY)\psi+\Tilde{\gamma}(Y)\nb^\Omega_X\psi.
\end{eqnarray}
The extrinsic Dirac operator on $M$ is then the first order elliptic differential operator acting on $\sm M$ defined by $\D^\Omega:=\mult^\Omega\circ\nb^\Omega$. 

\vspace{0.5cm}

Assume now that $\Omega^{n+1}$ has a smooth boundary $\Sigma^n:=\partial\Omega^{n+1}$. The orientation of $\Omega$ induces an orientation of $\Sigma$ which provides the existence of a unit vector field $N\in\Gamma(T\Omega_{|\Sigma})$ normal to $\Sigma$ and pointing inward $\Omega$. The existence of this vector field allows to induce on $\Sigma$ a spin structure from the one over $\Omega$. It follows that the bundle $\sm\Sigma:=\sm \Omega_{|\Sigma}$ is well-defined and is endowed with a spinorial Levi-Civita connection $\nb^\Sigma$, a Clifford multiplication $\mult^\Sigma$ and a Hermitian inner product $\<\,,\,\>$. The Clifford multiplication $\mult^\Sigma$ is defined for $X\in\Gamma(T\Sigma)$ and $\psi\in\Gamma(\sm\Sigma)$ by
\begin{eqnarray*}
\mult^\Sigma(X)\psi=\mult^\Omega(X)\mult^\Omega(N)\psi
\end{eqnarray*}
and is related, from (\ref{idmultclif}), to the space-time Clifford multiplication by 
\begin{eqnarray*}
\mult^\Sigma(X)\psi=\Tilde{\gamma}(X)\Tilde{\gamma}(N)\psi.
\end{eqnarray*}
Using this identification, the spin Gau{\ss} formula for the embedding $\Sigma^n\hookrightarrow\Omega^{n+1}$ reads as
\begin{eqnarray}\label{spin-gaus}
\nb^\Omega_X\psi=\nb^\Sigma_X\psi+\frac{1}{2}\Tilde{\gamma}\big(AX\big)\Tilde{\gamma}(N)\psi,
\end{eqnarray}
for all $\psi\in\Gamma(\sm\Sigma)$, $X\in\Gamma(T\Sigma)$ and where $A$ is the shape operator of the embedding of $\Sigma$ in $\Omega$ whose expression is given by (\ref{shape2}). The extrinsic Dirac operator of $\Sigma$ acting on $\sm\Sigma$ is defined by ${\D}^\Sigma:=\mult^\Sigma\circ\nb^\Sigma$ whose local expression is
\begin{eqnarray*}
\D^\Sigma\psi=\sum_{j=1}^n\mult^\Sigma(e_j)\nb^\Sigma_{e_j}\psi=\sum_{j=1}^n\Tilde{\gamma}(e_j)\Tilde{\gamma}(N)\nb^\Sigma_{e_j}\psi
\end{eqnarray*}
for all $\psi\in\Gamma\big(\sm\Sigma)$ and where $\{e_1,\cdots,e_n\}$ is a local orthonormal frame of $T\Sigma$. It is then straightforward to check that ${\D}^\Sigma$ is a first order elliptic linear differential operator which is formally self-adjoint for the $L^2$-scalar product on $\sm\Sigma$. On the other hand, a direct calculation using (\ref{compatibility1}) and the Gau\ss{} formula (\ref{spin-gaus}) gives the skew-commu\-ta\-ti\-vity rule
\begin{equation}\label{DNcommutes}
\D^\Sigma\big(\Tilde{\gamma}(N)\psi\big)=-\Tilde{\gamma}(N)\D^\Sigma\psi
\end{equation}
for any spinor field $\psi\in\Gamma({\sm}\Sigma)$.


\section{Spin outer untrapped submanifolds and eigenvalue estimates}



\subsection{The Dirac-Witten operator}\label{DW}


We first recall some standard facts about the Dirac-Witten operator first introduced by Witten \cite{Wi} in his proof of the positive energy theorem. Here we focus on the boundary expression of the associated Schr\"odinger-Lichnerowicz formula obtained by Gibbons, Hawking, Horowitz and Perry \cite{GHHP} in the context of positive energy theorems for black holes in the $(3+1)$ dimensional case. The Dirac-Witten operator is the first order elliptic differential operator acting on $\sm\Omega$ locally given by
\begin{eqnarray*}
\mathcal{D}\psi=\sum_{j=1}^{n+1}\Tilde{\gamma}(e_j)\widetilde{\nabla}_{e_j}\psi
\end{eqnarray*}
where $\{e_1,\cdots,e_{n+1}\}$ is a local orthonormal frame of the tangent bundle of $\Omega$. It is by now well-known that the following Schr\"odinger-Lichnerowicz type formula holds for $\psi\in\Gamma(\sm\Omega)$
\begin{eqnarray}\label{schrodinger}
\mathcal{D}^2\psi=\widetilde{\nabla}^*\widetilde{\nabla}\psi+\mathcal{R}\psi
\end{eqnarray}
where $\mathcal{R}$ is the endomorphism of $\sm\Omega$ defined by
\begin{eqnarray*}
\mathcal{R}\psi:=\frac{1}{4}\Big(R-|B|^2+\big({\rm Tr}(B)\big)^2-2\,\Tilde{\gamma}\big(\div\big(B-{\rm Tr}(B)I\big)\big)\Tilde{\gamma}(T)\Big)\psi
\end{eqnarray*}
and 
\begin{eqnarray}\label{adjointlaplacien}
\widetilde{\nabla}_{X}^*\psi=-\widetilde{\nabla}_{X}\psi+\Tilde{\gamma}(BX)\Tilde{\gamma}(T)\psi
\end{eqnarray}
is the $L^2$-formal adjoint of the spin connection $\widetilde{\nabla}$. On the other hand, using the compatibility property of the spin connection $\widetilde{\nabla}$ with the Hermitian product $(\;,\;)$ and the Stokes formula, we derive the following integration by parts formula
\begin{eqnarray}\label{ipp}
\int_{\Omega}\<\mathcal{D}\psi,\varphi\>d\Omega-\int_{\Omega}\<\psi,\mathcal{D}\varphi\>d\Omega=-\int_{\Sigma}\<\Tilde{\gamma}(N)\psi,\varphi\>d\Sigma
\end{eqnarray}
for all $\psi,\varphi\in\Gamma(\sm\Omega)$ and where $d\Omega$ (resp. $d\Sigma$) is the volume form of $\Omega$ (resp. $\Sigma$). Now using (\ref{adjointlaplacien}) and (\ref{ipp}) when integrating (\ref{schrodinger}) on $\Omega$ gives
\begin{eqnarray*}
\int_{\Omega}\Big(|\widetilde{\nabla}\psi|^2+\<\mathcal{R}\psi,\psi\>-|\mathcal{D}\psi|^2\Big)\,d\Omega=-\int_{\Sigma}\<\widetilde{\gamma}(N)\mathcal{D}\psi+\widetilde{\nabla}_N\psi,\psi\>\,d\Sigma.
\end{eqnarray*} 
A consequence of the dominant energy condition (\ref{DEC}) is the non-negativity of the endomorphism $\mathcal{R}$. Indeed, it is straightforward to check that the lowest eigenvalue of $\mathcal{R}$ is 
$$\mathcal{T}(T,T)-\sqrt{\sum_{j=1}^{n+1}\mathcal{T}(e_j,T)^2}=\mu-|J|,$$
which is non negative by taking $v=J$ in (\ref{DEC}). So if the dominant energy condition holds on $\E$, the previous formula reads
\begin{eqnarray*}
-\int_{\Omega}|\mathcal{D}\psi|^2\,d\Omega\leq-\int_{\Sigma}\<\Tilde{\gamma}(N)\mathcal{D}\psi+\widetilde{\nabla}_N\psi,\psi\>\,d\Sigma
\end{eqnarray*}
with equality if and only if $\psi$ satisfies $\widetilde{\nabla}_X\psi=0$ for all $X\in\Gamma(T\Omega)$. Following \cite{BaumMuller}, such a spinor field is referred to as an {\it extrinsic imaginary Codazzi spinor} (or $({\rm EIC})$-spinor) since using the Gau{\ss} formula (\ref{spingaussformule}) and the Clifford multiplications identification (\ref{idmultclif}), we observe that
\begin{equation}\label{EBK}\tag{${\rm EIC}$}
\widetilde{\nabla}_X\psi=0\quad\Longleftrightarrow\quad\nb^\Omega_X\psi=\frac{i}{2}\,\mult^\Omega(BX)\psi
\end{equation}
for all $X\in\Gamma(T\Omega)$.

On the other hand, if $\{e_1,\cdots,e_n,e_{n+1}=N\}$ is a local orthonormal frame of $T\Omega_{|\Sigma}$, we compute using the Gau{\ss} formula (\ref{spingaussformule}) and (\ref{spin-gaus}):
\begin{eqnarray*}
-\widetilde{\nabla}_{N}\psi-\Tilde{\gamma}(N)\mathcal{D}\psi=\EMD\psi-\frac{1}{2}H\psi-\frac{1}{2}\Tilde{\gamma}(BN)\Tilde{\gamma}(T)\psi.
\end{eqnarray*}
where 
\begin{eqnarray}\label{EMD}
\EMD\psi:=\D^\Sigma\psi+\frac{1}{2}K\Tilde{\gamma}(N)\Tilde{\gamma}(T)\psi
\end{eqnarray}
for all $\psi\in\Gamma(\sm\Sigma)$ and $BN\in\Gamma(T\Sigma)$ is defined by:
\begin{eqnarray*}
\<BN,X\>=B(X,N)
\end{eqnarray*}
for all $X\in\Gamma(T\Sigma)$. We have in particular, the following important formula
\begin{eqnarray}\label{FormEMD}
\EMD\psi=\frac{1}{2}H\psi+\frac{1}{2}\Tilde{\gamma}(BN)\Tilde{\gamma}(T)\psi+\sum_{j=1}^n\Tilde{\gamma}(e_j)\Tilde{\gamma}(N)\widetilde{\nabla}_{e_j}\psi
\end{eqnarray}
for all $\psi\in\Gamma(\sm\Sigma)$. We also note that since $\EMD$ is a zero-order modification of the extrinsic Dirac operator $\D^\Sigma$, it is obvious to check that it defines an elliptic linear differential operator of order one. Moreover from (\ref{MultRule}), one sees that the endomorphism $\Tilde{\gamma}(N)\Tilde{\gamma}(T)$ is pointwise self-adjoint with respect to the Hermitian scalar product $\<\;,\;\>$ so that $\EMD$ is symmetric for the $L^2$-scalar product on $\sm\Sigma$. 

\begin{remark}\label{Rem1}
{\rm 
The normal bundle $\mathcal{N}\Sigma$ of $\Sigma$ in $\E$ is a rank-two vector bundle over $\Sigma$ with induced metric of signature $(-,+)$. Moreover, the induced connection $\widetilde{\nabla}^\perp$ on $\mathcal{N}\Sigma$ is defined to be the fiberwise orthogonal projection of the space-time connection $\widetilde{\nabla}$ onto $\mathcal{N}\Sigma$. Fix any section $\nu\in\Gamma(\mathcal{N}\Sigma)$ that is outward-spacelike and of unit length. Then, there exists a unique $1$-form $\alpha_\nu$ on $\Sigma$ so that for all $X\in\Gamma(T\Sigma)$ we have
\begin{eqnarray*}
\alpha_\nu(X)=\<\widetilde{\nabla}^\perp_X\nu,\nu^\perp\>
\end{eqnarray*}
where the map $\nu\mapsto\nu^\perp$ is the involutive linear isomorphism of $\mathcal{N}\Sigma$ defined as follow: for $p\in\Sigma$ and $u,v\in\mathcal{N}_p\Sigma$ comprising an orthonormal basis, with $u$ future-directed and $v$ outward-spacelike, define $u^\perp=v$ and $v^\perp=u$ and extend linearly. The $1$-form $\alpha_\nu$ is the {\it connection $1$-form} of the normal bundle which, with $\nu$, completely determine $\widetilde{\nabla}^\perp$. The vector field $BN$ can be expressed in term of the vector field dual to the one-form $\alpha_{-N}$. Indeed, since $N$ is inward-spacelike, we have $(-N)^\perp=T$ and $T^\perp=-N$ so that
\begin{eqnarray*}
\alpha_{-N}(X)=-\<\widetilde{\nabla}^\perp_XN,T\>=-\<\widetilde{\nabla}_XN,T\>=B(X,N)=\<BN,X\>
\end{eqnarray*}
for all $X$ tangent to $\Sigma$. This implies that 
\begin{eqnarray*}
BN=\big(\alpha_{-N}\big)^\sharp
\end{eqnarray*}
where $(\;\;)^\sharp:T^\star\Sigma\rightarrow\ T\Sigma$ is the musical isomorphism between the co-tangent and the tangent bundles of $\Sigma$. Note that, in the following, since $N$ is always assumed to be inward-pointing, we let $\alpha_N:=-\alpha_{-N}$.
}
\end{remark}

\vspace{0.5cm}

As we will see, it is enough for our purpose to consider spin outer untrapped submanifolds with $\alpha_N=0$. From Remark \ref{Rem1} and the previous discussion, we obtain the following Reilly-type identity
\begin{proposition}\label{modlich}
Let $\Sigma^n:=\partial\Omega^{n+1}$ a spin outer untrapped submanifold in a time-oriented space-time $\E^{n+1,1}$ satisfying the Einstein equations (\ref{Einstein}) and the dominant energy condition (\ref{DEC}). If $\alpha_{N}=0$ then the inequality 
\begin{eqnarray}\label{mod-weit-ineq}
-\int_\Omega |\mathcal{D}\psi|^2\,d\Omega \leq\int_\Sigma\langle\EMD\psi-\frac{1}{2}H\psi,\psi\>\,d\Sigma
\end{eqnarray}
holds for all $\psi\in\Gamma(\sm\Omega)$. Moreover, equality occurs if and only if $\psi$ is a $({\rm EIC})$-spinor. 
\end{proposition}

\pf 
It only remains to prove that if $\psi$ is a $({\rm EIC})$-spinor on $\Omega$, then equality occurs in (\ref{mod-weit-ineq}). Indeed, from (\ref{EBK}), we immediately get $\mathcal{D}\psi=0$ so that the left-hand side of (\ref{mod-weit-ineq}) vanishes. On the other hand, using (\ref{EBK}) and the fact that $\alpha_N=0$ in (\ref{FormEMD}) finally lead to
\begin{eqnarray*}
\EMD\psi=\frac{1}{2}H\psi,
\end{eqnarray*}
which implies that the right-hand side also vanishes and so equality holds.
\qed


\subsection{An elliptic boundary problem for the Dirac-Witten operator}\label{EBP}


As before, $\Sigma^n$ is the boundary hypersurface of an $(n+1)$-dimensional Riemannian spin compact manifold $\Omega^{n+1}$ lying in a spacelike hypersurface of a time-oriented space-time $\E^{n+1,1}$. We define
two pointwise projections 
\begin{eqnarray*}
P_\pm:\sm\Sigma\longrightarrow \sm\Sigma
\end{eqnarray*}
on the induced Dirac bundle $\sm\Sigma$ over the submanifold $\Sigma$ by letting
\begin{equation}\label{defP}
P_\pm\psi=\frac{1}{2}\big(\psi\pm i\,\widetilde{\gamma}(N)\psi\big),\qquad\forall \psi\in\Gamma(\sm\Sigma).
\end{equation}
This boundary condition has been introduced in the seventies in order to propose a model for elementary particles \cite{CJJTW,CJJT,J}. It is a well known fact that the two orthogonal projections $P_\pm$ defined on the spin bundle $\sm\Sigma$ in (\ref{defP}), provide local elliptic boundary conditions for the Dirac-Witten operator $\mathcal{D}$ of $\Omega^{n+1}$ (see \cite{BB,BC}). The ellipticity of the boundary conditions $P_+$ and $P_-$ and that of the Dirac-Witten operator $\mathcal{D}$ on $\Omega$ allow to assert that the boundary problems associated to the realization $(\mathcal{D},P_\pm)$ on the domain $\Omega^{n+1}$ is of Fredholm type with smooth spinor fields as solutions, if they existed. Note also that the integration by parts formula (\ref{ipp}) shows that none of the conditions provided by $P_\pm$ makes $\mathcal{D}$ a formally self-adjoint operator. Instead, one can easily see that the boundary {\em realizations} $(\mathcal{D},P_+)$ and $(\mathcal{D},P_-)$ of $\mathcal{D}$ are adjoint to each other. We first deduce 
\begin{proposition}\label{boun-prob}
The following two types of inhomogeneous problems for the Dirac-Witten operator $\mathcal{D}$ 
of a $(n+1)$-dimensional Riemannian spin compact manifold with smooth boundary $\Sigma$ lying in a spacelike hypersurface of a time-oriented space-time $\E$ 
\begin{equation}\label{conf-loca-boun-cond}
\left\{
\begin{array}{lll}
\mathcal{D}\psi&=\Phi \qquad& \hbox{ {\rm on} } \Omega \\
P_\pm(\psi_{|\Sigma})&=0 \qquad& \hbox{ {\rm on} }\Sigma
\end{array}
\right. 
\end{equation}
have a unique smooth solution for any $\Phi\in\Gamma(\sm\Omega)$.
\end{proposition}

\noindent
\pf  
The two realizations of ${D}$ associated
with the two boundary conditions $P_\pm$ are the two unbounded operators
\beQ
\mathcal{D}_\pm:\hbox{\rm Dom}\,\mathcal{D}_\pm=\{\psi\in H^1(\sm\Omega)\,|\,
P_\pm(\psi_{|\Sigma})=0\}\longrightarrow L^2(\sm\Omega)
\eeQ
where $H^1(\sm\Omega)$ stands for the Sobolev space of $L^2$-spinors with weak $L^2$ covariant derivatives. 
From (\ref{ipp}), it follows that for the adjoint, one has $(\mathcal{D}_\pm)^* = \mathcal{D}_\mp$. Moreover, if $\psi\in \hbox{Dom}\,\mathcal{D}_\pm$ is a solution of the corresponding homogeneous problem, the ellipticity of both the Dirac-Witten operator $\mathcal{D}$ and the boundary condition $P_\pm=0$ imply that $\psi$ is smooth. On the other hand, taking $\varphi=i\psi$ in (\ref{ipp}) and recalling that the metric on $\sm\Omega$ is Hermitian, we have 
$$
0=2\int_ \Omega\langle\mathcal{D}\psi,i\psi\rangle\,d\Omega=\int_\Sigma\langle\psi,
i\,\Tilde{\gamma}(N)\psi\rangle\,d\Sigma=\mp\int_\Sigma|\psi|^2\,d\Sigma.
$$
Then one sees that the smooth $\mathcal{D}$-harmonic spinor $\psi$ on the compact manifold $\Omega$ has a vanishing trace $\psi_{|\Sigma}$ along the boundary hypersurface $\Sigma$. On the other hand, from the Gau{\ss} formula (\ref{spingaussformule}), it is straightforward to see that 
\begin{eqnarray*}
\mathcal{D}\psi=0\quad\Longleftrightarrow\quad \D^\Omega\psi=-\frac{i}{2}\big({\rm tr}\,B\big)\psi,
\end{eqnarray*}
where $\D^\Omega$ is the extrinsic Dirac operator of $\Omega$. As a consequence of this fact we get that, since $\psi$ is non-trivial and $\Omega$ is connected, its zeroes set has Hausdorff measure at most $n-1$ (see \cite{Ba1}). This contradicts the fact that $\psi$ vanishes identically on $\Sigma$ and so we conclude that $\psi$ vanishes on the whole of $\Omega$. Then 
$$
\ker\mathcal{D}_\pm=\{0\}\qquad \hbox{and}\qquad \hbox{\rm coker}\,
\mathcal{D}_\pm\cong\ker (\mathcal{D}_\pm)^*=\ker\mathcal{D}_\mp=\{0\}.
$$ 
Then the two realizations $\mathcal{D}_\pm$ are invertible operators, hence if $\Phi\in\Gamma(\sm\Omega)$ is a smooth spinor field on $\Omega$, there exists a unique solution $\psi\in H^1(\sm\Omega)$ of (\ref{conf-loca-boun-cond}). The classical regularity results imply the required smoothness of the solution $\psi$.
\qed

As a consequence of the previous proposition, we prove that the associated non-homogeneous boundary value problem has a unique solution.  
\begin{proposition}\label{boun-prob2}
On a $(n+1)$-dimensional Riemannian spin compact manifold $\Omega$ with smooth boundary $\Sigma$ lying in a spacelike hypersurface of a time-oriented space-time $\E$, the following boundary problem
$$ 
\left\{
\begin{array}{lll}
\mathcal{D}\Psi &=0 \qquad&\hbox{ {\rm on} } \Omega   \\
P_\pm\Psi_{|\Sigma} & = P_\pm\varphi \qquad&\hbox{ {\rm along }}\Sigma  
\end{array}
\right. 
$$
has a unique smooth solution $\Psi\in\Gamma(\sm\Omega)$ for all $\varphi\in\Gamma(\sm\Sigma)$.
\end{proposition}

\pf
Extend $\varphi$ to a spinor field $\widehat{\phi}\in\Gamma(\sm\Omega)$. From Proposition \ref{boun-prob}, there exists a unique 
$\widehat{\psi}\in\Gamma(\sm\Omega)$ smooth solution to the boundary problem
$$ 
\left\{
\begin{array}{lll}
\mathcal{D}\widehat{\psi} &=-\mathcal{D}\widehat{\phi} \qquad&\hbox{ {\rm on} } \Omega   \\
P_\pm\widehat{\psi}_{|\Sigma} & = 0 \qquad&\hbox{ {\rm on} } \Sigma  
\end{array}
\right. 
$$
and so $\Psi=\widehat{\psi}+\widehat{\phi}\in\Gamma(\sm\Omega)$ is the desired solution. 

\qed


\subsection{A generalized Reilly formula}


First, it is important for the following to study the behavior of the projections $P_\pm$ with respect to the Dirac-type operator $\EMD$:
\begin{lemma}\label{lemma1} 
For any smooth spinor field $\psi\in\Gamma(\sm\Sigma)$ we have 
\begin{equation}\label{antiDP}
\EMD P_\pm\psi=P_\mp\EMD\psi
\end{equation}
and
\begin{equation}\label{DPpm}
\int_\Sigma\langle\EMD\psi,\psi\rangle\,d\Sigma=2 \int_\Sigma
\langle\EMD P_+\psi,P_-\psi\rangle\,d\Sigma.
\end{equation}
\end{lemma}

\pf 
From (\ref{DNcommutes}), it is straightforward to check that 
\begin{eqnarray*}
\D^\Sigma\big(P_\pm\psi\big)=P_\mp\D^\Sigma\psi
\end{eqnarray*}
for all $\psi\in\Gamma(\sm\Sigma)$. Moreover since
\begin{eqnarray*}
\Tilde{\gamma}(N)\Tilde{\gamma}(T)P_\pm\psi=P_\mp\big(\Tilde{\gamma}(N)\Tilde{\gamma}(T)\psi\big)
\end{eqnarray*}
the relation (\ref{antiDP}) follows from the definition (\ref{EMD}) of $\EMD$. For the second point, we remark that the pointwise orthogonality of the decomposition $\psi=P_+\psi+P_-\psi$ combined with the previous relations give
\begin{eqnarray*}
\int_\Sigma\langle\EMD\psi,\psi\rangle\,d\Sigma=\int_\Sigma\langle P_+\EMD\psi,P_+\psi\rangle\,d\Sigma+\int_\Sigma\langle P_-\EMD \psi,P_-\psi\rangle\,d\Sigma.
\end{eqnarray*}
Then using (\ref{antiDP}) and the self-adjointness of $\EMD$ in the first term of the right-hand side of this equality prove (\ref{DPpm}).
\qed

The proof of the main result of this section (see Theorem \ref{main}) relies on a generalization of the classical Reilly formula for spinors (see Formula (\ref{GeneralizedReilly})). This approach has been developed by the first two authors in the spin Riemannian setting in \cite{HM} (see also \cite{HMR}). It should also be noticed that such a formula was also derived in \cite{MTX} for the standard Reilly formula on functions to study the critical points of the Wang-Yau quasi-local energy. So here we first prove
\begin{proposition}
Let $\Sigma^n$ be a spin outer untrapped submanifold of codimension-two in a spacetime ${\mathcal E}^{n+1,1}$ which satisfies the Einstein equations (\ref{Einstein}) and the dominant energy condition (\ref{DEC}). If $\alpha_N=0$, then for all $\psi\in\Gamma(\sm\Sigma)$, we have:
\begin{eqnarray}\label{P+}
\int_\Sigma\Big(\frac{1}{H}|\EMD P_+\psi|^2-\frac{H}{4}|P_+\psi|^2\Big)\,d\Sigma\geq 0.
\end{eqnarray}
Moreover, equality occurs if and only if there exists a $({\rm EIC})$-spinor $\Psi\in\Gamma(\sm\Omega)$ such that
$P_+\Psi=P_+\psi$.
\end{proposition}

\pf
For a spinor field $\psi\in\Gamma(\sm\Sigma)$ consider the following boundary value problem
$$
\left\{
\begin{array}{lll}
\mathcal{D}\Psi & =0 \qquad&\hbox{ {\rm on} } \Omega \\
P_+\Psi & =P_+\psi \qquad&\hbox{ {\rm along} }\Sigma 
\end{array}
\right.
$$ 
for the Dirac-Witten operator $\mathcal{D}$ and the pointwise boundary condition $P_+$. Existence and uniqueness of a smooth solution $\Psi\in \Gamma(\sm\Omega)$ for this boundary problem is ensured by Proposition \ref{boun-prob2}. On the other hand, since we assume that the spacetime ${\mathcal E}$ satisfies the dominant energy condition and $\alpha_N=0$, we can apply the inequality (\ref{mod-weit-ineq}) to $\Psi$ which gives 
\begin{eqnarray*}
0\leq\int_\Sigma\Big(\langle \EMD\Psi,\Psi\rangle-\frac{1}{2}H |\Psi|^2\rangle\Big)d\Sigma.
\end{eqnarray*}
Then from the relations (\ref{antiDP}), (\ref{DPpm}) and the fact that $P_+\Psi=P_+\psi$ we have
\begin{eqnarray}\label{DemiFin}
0\leq \int_\Sigma\Big(2\<\EMD(P_+\psi),P_-\Psi\>-\frac{1}{2}H|P_+\psi|^2-\frac{1}{2}H|P_-\Psi|^2\Big)d\Sigma.
\end{eqnarray}
On the other hand, since $\Sigma$ is an outer untrapped submanifold we have $H=\theta_+-\theta_->0$ and then
\begin{eqnarray*}
0\leq\Big|\sqrt{\frac{2}{H}}\EMD(P_+\psi)-\sqrt{\frac{H}{2}}P_-\Psi\Big|^2
\end{eqnarray*}
implies that
\begin{eqnarray*}
2\<\EMD(P_+\psi),P_-\Psi\>-\frac{1}{2}H|P_-\Psi|^2\leq\frac{2}{H}|\EMD(P_+\psi)|^2.
\end{eqnarray*}
Putting this estimate in (\ref{DemiFin}) gives the desired estimate. If equality is achieved, we have equality in (\ref{mod-weit-ineq}) so that $\Psi$ is a $({\rm EIC})$-spinor with $P_+\Psi=P_+\psi$. Conversely, assume that there exists such a spinor field $\Psi\in\Gamma(\sm\Omega)$. Since $\widetilde{\nabla}_X\Psi=0$ for all $X\in\Gamma(T\Sigma)$ and $\alpha_N=0$, we have from (\ref{FormEMD}) that $\EMD\Psi=\frac{1}{2}H\Psi$. Formula (\ref{antiDP}) in Lemma \ref{lemma1} implies that $\EMD(P_\pm\Psi)=\frac{1}{2}H P_\mp\Psi$ and since $P_+\Psi=P_+\psi$ along $\Sigma$ we have $\EMD(P_+\psi)=\frac{1}{2}HP_-\Psi$. This implies that
\begin{eqnarray}\label{RecEq}
0\leq\int_\Sigma\Big(\frac{1}{H}|\EMD P_+\psi|^2-\frac{H}{4}|P_+\psi|^2\Big)\,d\Sigma=\frac{1}{4}\int_\Sigma H\Big(|P_-\Psi|^2-|P_+\psi|^2\Big)d\Sigma.
\end{eqnarray}
On the other hand, we compute 
\begin{eqnarray*}
\frac{1}{4}\int_\Sigma H|P_+\psi|^2d\Sigma & = & \frac{1}{2}\int_\Sigma\<\EMD(P_-\Psi),P_+\psi\>d\Sigma\\
& = & \frac{1}{2}\int_\Sigma\<P_-\Psi,\EMD(P_+\psi)\>d\Sigma\\
& = & \frac{1}{4}\int_\Sigma H|P_-\Psi|^2d\Sigma
\end{eqnarray*}
which when used in (\ref{RecEq}) allows us to conclude. 
\qed

It turns out that inequality (\ref{P+}) also holds for the $P_-$-projection since there is an obvious symmetry between these two projections, that is
\begin{eqnarray}\label{P-}
\int_\Sigma\Big(\frac{1}{H}|\EMD P_-\psi|^2-\frac{H}{4}|P_-\psi|^2\Big)\,d\Sigma\geq 0
\end{eqnarray}
for all $\psi\in\Gamma(\sm\Sigma)$. Then summing (\ref{P+}) and (\ref{P-}) using the commutation properties (\ref{antiDP}) leads to following key inequality
\begin{proposition}\label{IneqInt}
Let $\Sigma^n$ be a spin outer untrapped submanifold of codimension-two in a spacetime ${\mathcal E}^{n+1,1}$ which satisfies the Einstein equations (\ref{Einstein}) and the dominant energy condition (\ref{DEC}). If $\alpha_N=0$, then for all $\psi\in\Gamma(\sm\Sigma)$, we have:
\begin{eqnarray}\label{GeneralizedReilly}
\int_\Sigma\Big(\frac{1}{H}|\EMD\psi|^2-\frac{H}{4}|\psi|^2\Big)\,d\Sigma\geq 0.
\end{eqnarray}
Moreover equality occurs if and only if there exists two $({\rm EIC})$-spinors $\Psi$, $\Phi\in\Gamma(\sm\Omega)$ with $P_+\Psi=P_+\psi$ and $P_-\Phi=P_-\psi$.
\end{proposition}


\subsection{A conformal eigenvalue estimate}


In this section, following the work \cite{HM} of the first two authors, we give an interpretation of Proposition \ref{IneqInt} in terms of the first eigenvalue of a Dirac-type operator for a metric conformally related to the initial one. This interpretation allows us to characterize the equality case in a more precise way. For this, recall that (see \cite{Hit,BHMM} for example) if $\<\,,\,\>_H=H^2\<\,,\,\>$ is a metric conformally related to $\<\,,\,\>$, there exists a bundle isometry between the two spinor bundles $\sm\Sigma$ and $\sm\overline{\Sigma}$ corresponding to the same spin structure and to the two conformal metrics. This identification will be denoted by
\begin{eqnarray*}
\psi\in\Gamma(\sm\Sigma)\mapsto \overline{\psi}\in\Gamma(\sm\overline{\Sigma}).
\end{eqnarray*}
Then the two Dirac operators $\D^\Sigma$ and $\overline{\D}^\Sigma$ acting respectively on $\sm \Sigma$ and $\sm\overline{\Sigma}$ are related for any spinor field $\psi\in\Gamma(\sm \Sigma)$ by
\begin{equation}\label{conf-boun-dira}
{\overline{\D}}^\Sigma\big( H^{-\frac{n-1}{2}}\psi\big)=H^{-\frac{n+1}{2}}\overline{{\D^\Sigma}\psi}. 
\end{equation}
Now we define the linear differential operator $\overline{\EMD}$ acting on $\sm\overline{\Sigma}$ by 
\begin{eqnarray}\label{EMDH}
\overline{\EMD}\,\overline{\varphi}:=\overline{\D}^\Sigma\overline{\varphi}+\frac{1}{2}H^{-1}K\,\mathcal{I}\overline{\varphi}
\end{eqnarray}
where $\mathcal{I}$ is the pointwise symmetric endomorphism given by
\begin{eqnarray*}
\mathcal{I} : \overline{\varphi}\in\Gamma(\sm\overline{\Sigma})\mapsto \overline{\Tilde{\gamma}(N)\Tilde{\gamma}(T)\varphi}\in\Gamma(\sm\overline{\Sigma}).
\end{eqnarray*}
Since $\overline{\EMD}$ is a zero order modification of the hypersurface Dirac operator $\overline{\D}^\Sigma$, it is clear that it defines an elliptic and $L^2$-self-adjoint operator of order one. Moreover, since $\Sigma$ is compact, its spectrum is an unbounded sequence of real numbers whose first eigenvalue is denoted by $\lambda_1(\overline{\EMD})$. Note that the spectrum of $\overline{\EMD}$ is symmetric with respect to zero since the endomorphism $\overline{\varphi}\mapsto\overline{\widetilde{\gamma}(N)\varphi}$ maps an eigenspinor for $\overline{\EMD}$ associated with $\lambda$ to an eigenspinor for $\overline{\EMD}$ associated with $-\lambda$. So without loss of generalities, we can assume that $\lambda_1(\overline{\EMD})$ is a non negative real number. On the other hand, using the conformal covariance relation (\ref{conf-boun-dira}), it is a straightforward computation to check that
\begin{eqnarray*}
\overline{\EMD}(H^{-\frac{n-1}{2}}\overline{\psi})=H^{-\frac{n+1}{2}}\overline{\EMD\psi}.
\end{eqnarray*}
With this relation, we observe that inequality (\ref{GeneralizedReilly}) now reads 
\begin{eqnarray*}
\int_{\Sigma}\Big(|\overline{\EMD}\,\overline{\psi}|^2-\frac{1}{4}|\overline{\psi}|^2\Big)\,d\overline{\Sigma}\geq 0
\end{eqnarray*}
for all $\psi=H^{-\frac{n-1}{2}}\varphi\in\Gamma(\sm\Sigma)$ and where $d\overline{\Sigma}=H^nd\Sigma$ represents the Riemannian measure with respect to the metric $\<\,,\,\>_H$. In other words, using the Rayleigh characterization of $\lambda_1(\overline{\EMD})^2$, we get
\begin{theorem}\label{main}
Let $\Sigma^n$ be a spin outer untrapped submanifold of codimension-two in a spacetime ${\mathcal E}^{n+1,1}$ which satisfies the Einstein equations (\ref{Einstein}) and the dominant energy condition (\ref{DEC}). If $\alpha_N=0$ then
\begin{eqnarray*}
\lambda_1(\overline{\EMD})^2\geq\frac{1}{4}.
\end{eqnarray*} 
Equality holds if and only if there exists a $({\rm EIC})$-spinor on $\Omega$. In this case, the  eigenspace corresponding to $\lambda_1(\overline{\EMD})=\frac{1}{2}$ consists of restrictions to $\Sigma$ of $({\rm EIC})$-spinor fields on $\Omega$ multiplied by the function $H^{-\frac{n-1}{2}}$. 
\end{theorem}

\pf 
We are left to prove the equality case. Assume now that equality is achieved, then there exists an eigenspinor $\overline{\psi}\in\Gamma(\sm\overline{\Sigma})$ for $\overline{\EMD}$ associated with eigenvalue $\frac{1}{2}$ which from (\ref{conf-boun-dira}) translates to $\EMD\varphi=\frac{1}{2}H\varphi$ with $\varphi=H^{\frac{n-1}{2}}\psi\in\Gamma(\sm\Sigma)$. Therefore we are in the equality case of Proposition \ref{IneqInt} and so there exists two $({\rm EIC})$-spinors  $\Psi$, $\Phi\in\Gamma(\sm\Sigma)$ with
$P_+\Psi=P_+\varphi$, $P_-\Phi=P_-\varphi$. In particular, from the equation (\ref{EBK}), the fact that $\alpha_N=0$ and the relation (\ref{FormEMD}), we obtain
\begin{eqnarray*}
\EMD\Psi=\frac{1}{2}H\Psi\quad\text{and}\quad\EMD\Phi=\frac{1}{2}H\Phi
\end{eqnarray*}
Now from Lemma \ref{lemma1} and the two previous equations, we easily deduce that
\begin{eqnarray*}
\EMD(P_\pm\Psi)=\frac{1}{2}H P_\mp\Psi\quad\text{and}\quad\EMD(P_\pm\Phi)=\frac{1}{2}H P_\mp\Phi
\end{eqnarray*}
so that
\begin{eqnarray*}
\frac{1}{2}H\,P_-\Psi=\EMD(P_+\Psi)=\EMD(P_+\varphi)=P_-(\EMD\varphi)=\frac{1}{2}H\,P_-\varphi=\frac{1}{2}H\,P_-\Phi.
\end{eqnarray*}
This implies that there exists a $({\rm EIC})$-spinor $\Psi\in\Gamma(\sm\Omega)$ such that $\Psi_{|\Sigma}=\varphi$. The converse is clear.
\qed


\section{CNNC codimension-two submanifolds}\label{AlexMink1}


In this section, we will see that the proof of Theorem \ref{AlexMink} is a a direct consequence of the conformal eigenvalue estimate proved in Theorem \ref{main}. Recall that a codimension-two submanifold of a Lorentzian manifold has CNNC if there exists a future null normal vector field $\Tilde{\mathcal{L}}$ such that $\Sigma$ is torsion-free with respect to $\TL$ and $\<\mathcal{H},\Tilde{\mathcal{L}}\>$ is a constant. We can assume that $\TL:=\TL_-$ is an inner vector field since the case where $\Tilde{\mathcal{L}}$ is an outer one can be treated in a same way. Then the strategy of the proof is as follows: we start from $\Sigma$ which bounds a compact domain $\Omega$ of a spacelike hypersurface in $\mathbb{R}^{n+1,1}$ with $T$ and $N$ the unit timelike future-directed normal to $\Omega$ in $\mathbb{R}^{n+1,1}$ and the unit inner normal to $\Sigma$ in $\Omega$. We first deform the domain $\Omega$ in a another domain $\Tilde{\Omega}$ with $\Sigma$ as boundary and whose inner null expansion is precisely given by $\widetilde{\theta}_-=-\<\mathcal{H},\TL_-\>$. Then from the torsion-free assumption, we deduce that $\alpha_{\Tilde{N}}=0$ so that Theorem \ref{main} may apply in a subtle way to this situation.


\subsection{Deformations of a spacelike domain spanned by a codimension-two  untrapped submanifold}


Suppose that $\Sigma^n$ is an orientable co\-di\-men\-sion-two spacelike submanifold in a Lorentzian manifold $\E^{n+1,1}$ and consider its mean curvature vector field ${\mathcal H}$, which is a normal field on $\Sigma^n$. If the immersion of $\Sigma^n$ into $\E^{n+1,1}$ factorizes through a spacelike hypersurface $M^{n+1}$ of the spacetime, we know that this factorization yields corresponding null expansions $(\theta_+,\theta_-)$ of ${\mathcal H}$, where the smooth functions $\theta_\pm$ are linked by the relation $|{\mathcal H}|^2=-\theta_+\theta_-$ and are respectively (up to a constant) the components of ${\mathcal H}$ with respect to the pair of  null vectors fields $(T\pm N)$, with $T$ the unit timelike future-directed normal field to $M^{n+1}$ and $N$ a choice of unit spacelike normal field of $\Sigma^n$ in $M^{n+1}$. When $\Sigma^n$ bounds a domain $\Omega^{n+1}$ in $M^{n+1}$, that is, $\Sigma^n=\partial\Omega^{n+1}$, we choose $N$ as the inner normal. 

Assume now that the same codimension-two spacelike submanifold $\Sigma^n$ factorizes through another spacelike hypersurface $P^{n+1}$ of $\E^{n+1,1}$. The new factorization provides us a different Lorentzian orthonormal reference $\{\Tilde{T},\Tilde{N}\}$ for the normal plane of $\Sigma$ in $\E$. In fact, it is obvious that there must be a function $f\in C^\infty(\Sigma)$, such that
\begin{equation}\label{LorentzTrans}
\Tilde{T}=(\cosh f) T-(\sinh f) N,\qquad \Tilde{N}=-(\sinh f) T+(\cosh f) N.
\end{equation}
The new Lorentzian reference also determines two null normal vectors $(\Tilde{T}\pm\Tilde{N})$ and the corresponding new  null expansions of ${\mathcal H}$ are given by 
$$
\Tilde{\theta}_+=e^f\theta_+,\qquad \Tilde{\theta}_-=e^{-f}\theta_-,
$$
or equivalently
\begin{equation}\label{thetachange}
\Tilde{K}+\Tilde{H}=e^f(K+H),\qquad \Tilde{K}-\Tilde{H}=e^{-f}(K-H).
\end{equation}
Conversely, it is important to remark that, provided that $\Sigma$ spans a compact spacelike domain $\Omega$ of $\E$, any choice of $f\in C^\infty (\Sigma)$ determines a factorization through another spacelike hypersurface which can be obtained from $\Omega$ through a suitable deformation. We refer to \cite{HMR2} for a proof of this result.

\begin{lemma}\label{theorem0} 
Let $\Sigma$ be a compact spacelike codimension-two submanifold embedded in the Lorentzian manifold $\E$ and ${\mathcal H}$ its mean curvature vector field. Suppose that $\Sigma$ spans a compact spacelike domain $\Omega$ in $\E$ and that $(\theta_+$, $\theta_-)$ are the outer and inner null expansions of ${\mathcal H}$ corresponding to the embedding of $\Sigma$ in $\Omega$. For any smooth function $f\in C^\infty(\Sigma)$, there exists a compact spacelike hypersurface $\Omega_f$ in $\E$ such that $\partial\Omega_f=\Sigma$ and the corresponding null expansions of ${\mathcal H}$ are given by
$$
\theta^f_+=e^f\theta_+,\qquad \theta^f_-=e^{-f}\theta_-.
$$
\end{lemma}

\begin{remark}\label{remark2}{\rm 
Suppose that the spacelike hypersurface $\Omega^{n+1}$ spanned by the codimension-two spacelike submanifold $\Sigma$ is endowed with a spin structure (this always occurs when the spacetime itself is spin). Since the new hypersurface $\Omega_f$ built in Lemma \ref{theorem0} above for a given $f\in C^\infty (\Sigma)$ is obtained by slightly deforming $\Omega$ near its boundary, all these $\Omega_f$ are homotopic to each other, and so they induce the same spin structure on $\Sigma$. It is obvious that all their Riemannian metrics, which come from the Lorentzian metric on ${\mathcal E}$, determine the same Riemannian metric on $\Sigma$.} 
\end{remark}

Now a particular choice of the function $f\in C^\infty (\Sigma)$ will allow us to factorize this embedding through a spacelike hypersurface $\widetilde{\Omega}$ of $\E$ specially adapted to obtain the proof of Theorem \ref{AlexMink}.
 
\begin{proposition}\label{goodK} 
Let $\Sigma$ be an outer untrapped spacelike codimension-two submanifold embedded in a Lorentzian manifold $\E$. Then, there exists a smooth function $f\in C^\infty(\Sigma)$ such that the corresponding $\Tilde{\Omega}:=\Omega_f$ given by Lemma \ref{theorem0} above has null expansions of ${\mathcal H}$ given by 
$$
\widetilde{\theta}_+=\frac{|\mathcal{H}|^2}{\<\mathcal{H},\TL_-\>}\qquad\widetilde{\theta}_-=-\<\mathcal{H},\TL_-\>
$$
where $\TL_-$ is an inner future-directed null normal. Moreover, if $\Sigma$ is torsion-free with respect to $\TL_-$ and if $\{\Tilde{T},\Tilde{N}\}$ denotes the associated Lorentzian frame obtained by the transformation (\ref{LorentzTrans}),  we have 
$$
\alpha_{\Tilde{N}}(X)=\<\widetilde{\nabla}_X\widetilde{N},\widetilde{T}\>=0
$$  
for all $X\in\Gamma(T\Sigma)$. 
\end{proposition}

\pf
Let $(\theta_+$, $\theta_-)$ be the outer and inner null expansions of the mean curvature vector field ${\mathcal H}$ of $\Sigma$ corresponding to the embedding of $\Sigma$ in $\Omega$ and to a given orientation of $\Omega$. Since $\Sigma$ is an outer untrapped submanifold we have that $-\theta_+\theta_-=|{\mathcal H}|^2>0$ and $\theta_+-\theta_-=2H>0$. So $\theta_+$ and $-\theta_-$ are positive everywhere on $\Sigma$. Then, we can put 
$$
f=\log\Big(-\frac{\theta_-}{\<\mathcal{H},\TL_-\>}\Big)
$$
in Lemma \ref{theorem0} and it suffices to define $\widetilde{\Omega}=\Omega_f$ for this choice of function $f$. On the other hand, it is obvious to observe that in this case, we have
\begin{eqnarray*}
\TL_-=\widetilde{T}+\widetilde{N}
\end{eqnarray*}
and a straightforward computation proves that
\begin{eqnarray*}
(\widetilde{\nabla}_X\TL_-)^\perp=\widetilde{B}(X,\widetilde{N})\,\TL_-
\end{eqnarray*}
for all $X\in\Gamma(T\Sigma)$. Now since $\Sigma$ is torsion-free with respect to $\TL_-$, we have $(\widetilde{\nabla}\TL_-)^\perp=0$ so that $\widetilde{B}(X,\widetilde{N})=0$ for all $X\in\Gamma(T\Sigma)$ which, as noticed in Remark \ref{Rem1}, is equivalent to the fact that $\alpha_{\Tilde{N}}(X)=0$ for all $X\in\Gamma(T\Sigma)$, as claimed.
\qed


\subsection{CNNC codimension-two submanifolds in the Minkowski spacetime}


In this section, we give the proof of Theorem \ref{AlexMink}. We first remark that compact untrapped submanifolds, that is, compact codimension-two spacelike submanifolds of Minkowski spacetime having spacelike mean curvature vector field, must be mean-convex in any spacelike domain that they might bound. Indeed we have
\begin{lemma}\label{2classes}
Let $\Sigma^n$ be a compact untrapped codimension-two submanifold in the Minkowski spacetime ${\mathbb R}^{n+1,1}$. Suppose that $\Sigma^n$ is the boundary of a  spacelike  domain in ${\mathbb R}^{n+1,1}$. Then $\Sigma^n$ is  mean-convex in
this domain, that is, $\Sigma$ is an outer untrapped submanifold.
\end{lemma}

\pf
Let $\Omega$ be a spacelike hypersurface  in $\R^{n+1,1}$ and spanned by $\Sigma$. Suppose that $\theta_+=H+K$ and 
$\theta_-=K-H$ are the corresponding null expansions of the mean curvature vector ${\mathcal H}$, after choosing suitable $T$ and $N$. Choose a unit timelike vector $a\in\R^{n+1,1}$ and project $\Omega^{n+1}$ onto the Euclidean slice $\R^{n+1}\subset
\R^{n+1,1}$ orthogonal to $a$. Since $\Omega$ is spacelike, we obtain an immersion $\phi$ from $\Omega$ into $\R^{n+1}$. It is not so difficult to check that the mean curvature $H'$ of the immersed hypersurface given by $\phi_{|\Sigma}:\Sigma\rightarrow {\mathbb R}^{n+1}$ with respect to the inner orientation is given by
$$
H'=H\cosh f +K\sinh f=e^f\theta_+-e^{-f}\theta_-,
$$
where $f\in C^\infty (\Sigma)$ is the smooth function given by
$$
f=\arg\cosh\frac{\langle T,a\rangle}{\sqrt{\langle T,a\rangle^2-\langle N,a\rangle^2}}.
$$ 
As a first consequence, $H'^2\ge |{\mathcal H}|^2$. As we are supposing that $\Sigma$ has spacelike mean curvature, we
conclude that $H'$ does not change its sign on each component of $\Sigma^n$. A second consequence is that $\Sigma$ is mean-convex, that is, $H>0$,  if and only if $H'>0$. But this is precisely the case because each component of $\Sigma$ is compact and any compact hypersurface in a Euclidean space has at least an elliptic point. Thus, $\Sigma$ must be mean-convex and $H'\ge |{\mathcal H}|>0$.
\qed

\noindent {\textit {Proof  of {\bf Theorem \ref{AlexMink}} :}}
Let $\Sigma^n$ be an untrapped CNNC submanifold of codimension-two in $\R^{n+1,1}$ which bounds a compact spacelike domain $\Omega^{n+1}$ of $\R^{n+1,1}$. From the previous lemma and the fact that $\R^{n+1,1}$ is spin, we deduce immediately that $\Sigma$ is in fact a spin outer untrapped submanifold. Since $\Sigma$ has CNNC, it is a torsion-free submanifold with respect to a future null normal vector field $\TL_-$. Once again, we will assume that $\TL_-$ is inward-pointing, the outward-pointing case can be treated in a same way. Now using Proposition \ref{goodK}, there exists a compact spin spacelike domain $\widetilde{\Omega}^{n+1}$ with smooth boundary $\Sigma$ such that $\widetilde{\theta}_-=-\<\mathcal{H},\TL_-\><0$. Moreover, if $\widetilde{N}$ (resp. $\widetilde{H}$) denotes the inner unit vector field normal to (resp. the mean curvature of) $\Sigma$ in $\widetilde{\Omega}$, we have, again applying Proposition \ref{goodK} (resp. Lemma \ref{2classes}), that $\alpha_{\widetilde{N}}=0$ (resp. $\widetilde{H}>0$). So if $\overline{\EMD}$ denotes the Dirac-type operator defined in (\ref{EMDH}) for the metric $\<\,,\,\>_{\widetilde{H}}=\widetilde{H}^2\<\,,\,\>$, we can apply Theorem \ref{main} to deduce that $\lambda_1(\overline{\EMD})\geq 1/2$. On the other hand, since $\R^{n+1,1}$ carries parallel spinor fields it induces $({\rm EIC})$-spinor fields on $\widetilde{\Omega}$ and so it is immediate, again from Theorem \ref{main}, that we have in fact $\lambda_1(\overline{\EMD})=1/2$. Consider now a parallel spinor field $\Phi\in\Gamma(\mathbb{S}\R^{n+1,1})$ and define 
\begin{eqnarray*}
\Psi:=\widetilde{\gamma}\Big(\frac{1}{n}\widetilde{\theta}_-\xi-\TL_-\Big)\Phi\in\Gamma(\sm\Sigma)
\end{eqnarray*}
where $\xi$ denotes the position vector field in $\R^{n+1,1}$. Now since $\Sigma$ has CNNC we have that $\<\mathcal{H},\TL_-\>$
is constant and from our choice of $\widetilde{\Omega}$, this implies that $\widetilde{\theta}_-$ is a (negative) constant. From this fact, we compute using (\ref{FormEMD}) that 
\begin{eqnarray*}
\EMD\Psi=\frac{1}{2}\widetilde{H}\Psi
\end{eqnarray*}
and so it induces an eigenspinor for $\overline{\EMD}$ associated with the eigenvalue $(1/2)$. From the equality case of Theorem \ref{main} we deduce that $\Psi$ is the restriction to $\Sigma$ of an $({\rm EIC})$-spinor on $\widetilde{\Omega}$. Then, for all $X\in\Gamma(T\Sigma)$, we have
\begin{eqnarray*}
0=\widetilde{\nabla}_X\Psi=\widetilde{\gamma}\Big(\frac{1}{n}\widetilde{\theta}_-X-\widetilde{\nabla}_X\TL_-\Big)\Phi
\end{eqnarray*}
and since $\Phi$ has no zero, we immediately get that the null second fundamental form with respect to $\TL_-$, defined by (\ref{NSFF}), satisfies
\begin{eqnarray*}
\chi_-(X,Y)=\frac{1}{n}\widetilde{\theta}_-\<X,Y\>
\end{eqnarray*}
for all $X$, $Y\in\Gamma(T\Sigma)$, that is $\Sigma$ lies in a shearfree null hypersurface.
\qed



\begin{thebibliography}{BHHM}

\bibitem[A]{A}
A.D. Alexandrov, {\em A characteristic property of spheres}, Ann. Mat. Pure Appl, {\bf 58} (1962), 303-315.

\bibitem [Ba]{TB-Baum}
H. Baum, {\em Spin-Strukturen und Dirac-Operatoren {\"u}ber pseudo-Riemann\-sche Mannig\-faltig\-keiten}, Teubner-Verlag, Leipzig, 1981.

\bibitem [BM]{BaumMuller}
H. Baum, O. M\"uller, {\em Codazzi spinors and globally hyperbolic manifolds with special holonomy}, 
Math. Zeit, {\bf 258} (1) (2008), 185--211.

\bibitem [B{\"a}1]{Ba1} 
C. B{\"a}r, {\em On nodal sets for Dirac and Laplace operators}, 
Commun. Math. Phys., {\bf 188} (1997), 709--721.  

\bibitem [B{\"a}Ba]{BB} 
C. B{\"a}r, W. Ballmann, {\em Boundary value problems for elliptic differential operators of first order}, Surveys in Differential Geometry {\bf 17} (2012), 1--78.

\bibitem [BC]{BC} 
R. Bartnik, P. Chru{\'s}ciel, {\em Boundary value problems
for Dirac-type equations}, J. reine angew. Math., {\bf 579} (2005), 13--73. 

\bibitem [BHMM]{BHMM} J.P. Bourguignon, O. Hijazi,
J.-L. Milhorat, A. Moroianu, {\it A Spinorial Approach to
Riemannian and Conformal Geometry}, Monographs in Mathematics, EMS, Z\"urich, 2015.

\bibitem [BJ1]{BJ1}
H.L. Bray, J.L. Jaurequi, {\em Time flat surfaces and the monotonicity of the spacetime Hawking mass}, Comm. Math. Phys., {\bf 335} (1) (2015), 285--307.

\bibitem [BJ2]{BJ2}
H.L. Bray, J.L. Jaurequi, {\em On curves with nonnegative torsion}, Arch. Math., {\bf 104} (6) (2015), 561--575.
 
\bibitem[Br]{Br}
S. Brendle, {\em Constant mean curvature surfaces in warped product manifolds}, Publications Math\'ematiques de l'IHES {\bf 117} (2013), 247--269.

\bibitem [CWW]{CWW}
P.N. Chen, M.T. Wang and Y.K. Wang, {\em Rigidity of time-flat surfaces in the Minkowski spacetime}, Math. Res. Lett, {\bf 21} (6) (2014), 1227--1240.

\bibitem [CJJTW]{CJJTW} A. Chodos, R.L. Jaffe, K. Johnson, C.B. Thorn, V.F. Weisskopf, {\em New extended model of hadrons}, Phys. Rev. D, {\bf 9} (1974), 3471--3495.

\bibitem [CJJT]{CJJT} A. Chodos, R.L. Jaffe, K. Johnson, C.B. Thorn, {\em Baryon structure in the bag theory}, Phys. Rev. D, {\bf 10} (1974), 2599--2604.

\bibitem [GHHP]{GHHP}
G. Gibbons, S. Hawking, G. Horowitz and M. Perry, {\em Positive mass theorems for black holes}, Comm. Math. Phys., {\bf 88} (1983), 295--308.

\bibitem [HM]{HM} 
O. Hijazi, S. Montiel, {\em A holographic principle for the existence
of parallel spinor fields and an inequality of Shi-Tam type}, Asian J. Math. {\bf 18} (3) (2014), 489-506.

\bibitem [HMR1]{HMR}
O. Hijazi, S. Montiel, S. Raulot, {\em On a Liu-Yau type inequality for surfaces}, Pac. J. Math., {\bf 272} (1) (2014), 177-199.

\bibitem [HMR2]{HMR2}
O. Hijazi, S. Montiel, S. Raulot, {\em Dirac operators on time flat submanifolds with applications}, Comm. Math. Phys., {\bf 351} (3) (2017), 1177-1194.

\bibitem [HMRo]{HMRo}
O. Hijazi, S. Montiel, A. Rold{\'a}n, {\em Dirac operators on hypersurfaces of manifolds
with negative scalar curvature}, Ann. Global Anal. Geom., {\bf 23} (2003), 247-264.

\bibitem [HMZ]{HMZ} 
O. Hijazi, S. Montiel, X. Zhang, {\em Dirac operator on embedded hypersurfaces}, 
Math. Res. Lett, {\bf 8} (2001), 195--208.

\bibitem [Hit]{Hit} N. Hitchin, {\em Harmonic spinors}, Adv. 
in Math., {\bf 14} (1974), 1--55.

\bibitem [J]{J} K. Johnson, {\em The M.I.T. bag model}, Acta Phys. Pol., {\bf B6} (1975), 
865--892.

\bibitem [L]{leistner}
T. Leistner, {\em On the classification of Lorentzian holonomy groups}, J. Diff. Geom., {\bf 76}, no.3 (2007), 423--484.
	
\bibitem [MTX]{MTX}
P. Miao, L.-F. Tam and N. Xie, {\em Critical points of the Wang-Yau quasi-local energy}, Ann. Henri Poincar\'e, {\bf 12} (2011), no. 5, 987--1017. 

\bibitem[Mo]{Mo}
S. Montiel, {\em Unicity of constant mean curvature hypersurfaces in some Riemannian manifolds}, Indiana Univ. Math. J., {\bf 48} (1999), 711--748.

\bibitem[MR]{MR}
S. Montiel, A. Ros {\em Compact hypersurfaces: the Alexandrov theorem for higher order mean curvatures}, Pitman Monographs and Surveys in Pure and Applied Mathematics {\bf 52} (1991) (in honor of M.P. do Carmo; edited by B. Lawson and K. Tenenblat), 279-296.

\bibitem [WWZ]{WWZ}
M.-T. Wang, Y.-K. Wang, X. Zhang, {\em Minkowski formulae and Alexandrov theorems in spacetime}, J. Differential Geom., {\bf 105} (2) (2017), 249--290.

\bibitem[Wi]{Wi} E. Witten, {\em A new proof of the positive energy theorem}, Commun. Math. Phys., {\bf
80} (1981), 381--402.
  
\end{thebibliography}
\end{document}